\documentclass[12pt]{amsart}

\theoremstyle{plain}
\newtheorem{theorem}{Theorem}[section]
\newtheorem{lemma}[theorem]{Lemma}
\newtheorem{proposition}[theorem]{Proposition}
\newtheorem{corollary}[theorem]{Corollary}

\theoremstyle{definition}
\newtheorem{definition}{Definition}

\theoremstyle{remark}
\newtheorem{remark}{Remark}

\def \Cset {{\mathbb C}}

\def \calO  {{\mathcal{O}}}

\def \calM {{\mathcal{M}}}
\def \calJ {{\mathcal{J}}}

\def \ra  {\rightarrow}           %maps

\def \la {\langle}
\def \ra {\rangle}
\def \ol {\overline}            

\def \rank { {\mathrm{rank}} }

\def \g  {\mathfrak{g}}   % Lie algebra letters

\renewcommand \max { {\mathrm{max}} }
\def \hs {\hspace{.2in}}
\def \sc {{\scriptscriptstyle C}}
\def \scd {{\scriptscriptstyle C_\delta}}
\def \Wc {W_\sc}
\def \Wcd {W_{\scd}}
\def \Wcdm {W^-_{\scd}}
\def \Wcm {W_{\sc}^-}
\def \mc {m_{\sc}}
\def \mcd {m_{\scd}}

\def \mcp {m_{{\scriptscriptstyle C'}}}

\def \nem {\neq \emptyset}

\def \sJm {{\scriptscriptstyle J_m}}
\def \sJmp {{\scriptscriptstyle J_m'}}

\def \sJ {{\scriptscriptstyle J}}
\def \mJ {m_{\sJ}}
\def \sK {{\scriptscriptstyle K}}
\def \mK {m_\sK}

\def \SL {SL(n+1, {\bf k})}

\input xy
\xyoption{all}

\begin{document}
%\tableofcontents
\setlength{\baselineskip}{1.2\baselineskip}
%%%%%%%%%%%%%%%%%%%%%%%%%%%%%%%%%%%%%%%%%%%%%%%%%%%%%%%%%%%%%%%%%%%%%%%%%%%
%%%%%%%%%%%%%%%%%%%%%%    Title    %%%%%%%%%%%%%%%%%%%%%%%%%%%%%%%%%%%%%%%%
\title[Intersections of conjugacy classes and Bruhat cells]
{On intersections of conjugacy classes and Bruhat cells}
\author[Kei Yuen Chan, Jiang-Hua Lu, Simon Kai-Ming To]{Kei Yuen Chan, Jiang-Hua Lu, Simon Kai-Ming To}
\address{
Department of Mathematics \\
Hong Kong University \\
Pokfulam Rd., Hong Kong}
\email{\!\!keiyuen@graduate.hku.hk,\!\! jhlu@maths.hku.hk,\!\! h0389481@graduate.hku.hk}

\begin{abstract} For a connected semi-simple algebraic group $G$ 
over an algebraically closed field ${\bf k}$ 
and a fixed pair $(B, B^-)$ 
of opposite Borel subgroups of $G$, 
we determine when 
the intersection of a conjugacy class $C$ in $G$ and a double
coset $BwB^-$ is nonempty, where $w$ is in the Weyl group $W$ of $G$. 
The question comes from 
Poisson geometry, and our answer 
is in terms of the Bruhat order on $W$ and an involution $\mc \in W$ associated to $C$.
We prove that  the element $\mc$ is the unique maximal length element in
its conjugacy class in $W$, and we classify all such elements in $W$.
For $G = SL(n+1, {\bf k})$, we 
describe $\mc$ explicitly for every conjugacy class $C$, and when $w \in W \cong S_{n+1}$ is an involution, we give an
explicit answer to when 
$C \cap (BwB)$ is nonempty.
\end{abstract}

\maketitle

%%%%%%%%%%%%%%%%%%%%   Introduction   %%%%%%%%%%%%%%%%%%%%%%%%%%%%%%%%%%%%%%%%
\section{Introduction}\label{sec-intro}
\subsection{The set up and the results}\label{subsec-1.1}

Let $G$ be a connected semi-simple algebraic group over an algebraically closed field 
${\bf k}$, and let
$B$ and $B^-$ be a pair of opposite Borel
subgroups of $G$. Then $H = B \cap B^-$ is a Cartan subgroup of $G$.
Let $W= N_G(H)/H$ be the Weyl group, where $N_G(H)$ is the 
normalizer of $H$ in $G$. One then has the well-known 
Bruhat decompositions
\[
G = \bigsqcup_{w \in W} BwB = \bigsqcup_{w \in W} BwB^- 
\hs (\mbox{disjoint unions}).
\]
Subsets of $G$ of the form $BwB$ or $Bw'B^-$, where $w, w' \in W$, will be called Bruhat cells in $G$.
For a subset $X$ of $G$, let $\overline{X}$ be the Zariski closure of $X$ in $G$.
The Bruhat order on $W$ is the partial order defined by
\[
w_1 \le w_2 \hs \Longleftrightarrow \hs Bw_1B \subset \ol{Bw_2B}, \hs w_1, w_2 \in W.
\]
Given a conjugacy class $C$ of $G$, let
\begin{eqnarray}\label{eq-Wc}
\Wc &=& \{w \in W: \; C\cap (BwB) \neq \emptyset\},\\
\label{eq-Wcm}
\Wcm &=& \{w \in W: \; C\cap (BwB^-) \neq \emptyset\}.
\end{eqnarray}
The sets $\Wc$ have been studied by several authors (see, for example, \cite{EG1, EG2} by E. Ellers and N. Gordeev and \cite{Ca1} by G. Carnovale)
and are not easy to determine even for the case of $G = SL(n, {\bf k})$ (see \cite{EG2}). 
On the other hand, for a conjugacy class $C$ in $G$, let $\mc$ be the unique element in $W$
such that $C \cap (B\mc B)$ is dense in $C$ (see Lemma~\ref{le-def-mc}). 
Our first result, Theorem~\ref{th-Wcm-mc}, states that,
for every conjugacy class $C$ in $G$,
\[
\Wcm = \{w \in W: \, w \leq \mc\}.
\]
Thus the set $\Wcm$ is completely determined by $\mc$ and the Bruhat order on $W$.

Theorem~\ref{th-Wcm-mc} is motivated by Poisson geometry. It is shown in 
\cite{EL} that when ${\bf k} = \Cset$,
the connected complex semi-simple Lie group $G$ carries a holomorphic Poisson structure $\pi_0$,
invariant under conjugation by elements of $H$, such that the
nonempty intersections $C \cap (B wB^-)$ are exactly 
the $H$-orbits of symplectic leaves of $\pi_0$, 
where $C$ is a conjugacy class
in $G$ and $w \in W$. To describe  the symplectic leaves of
$\pi_0$, one first needs to determine when 
$C \cap (BwB^-)$ is nonempty.
By \cite[Theorem 1.4]{R2}, 
the nonempty intersections $C \cap (BwB^-)$ are always  smooth and irreducible. 
The geometry of such intersections and applications to Poisson geometry will be carried out
elsewhere.

%In \cite{EG1999}, E. Ellers and N. Gordeev proved that for every non-central conjugacy class $C$ in $G$, the set $\Wcm$ contains
%all generalized Coxeter elements (see \cite[Definition 2.2]{EG1} or $\S$\ref{subsec-EG1999} for the definition).
%Discussions on this result of Ellers and Gordeev in light of our Theorem~\ref{th-Wcm-mc} are given in $\S$\ref{subsec-EG1999}.

The elements $\mc$ play an important role in the study of spherical conjugacy classes, i.e., conjugacy
classes in $G$ on which the $B$-action by conjugation has a dense orbit.  Let ${\bf k} = \Cset$ and let 
$\g$ be the Lie algebra of $G$. 
In connection with their
proof of the de Concini-Kac-Procesi conjecture on representations of the quantized universal enveloping algebra 
${\mathcal U}_{\epsilon}(\mathfrak g)$ at roots of unity over spherical conjugacy classes, 
N. Cantarini, G. Carnovale, and M. Costantini proved  \cite[Theorem 25]{CCC} that a conjugacy class $C$ in $G$ is spherical if and only if
$\dim C = l(\mc) + \rank(1-\mc)$,
 where $l$ is the length function on $W$, and
$\rank(1-\mc)$ is the rank of the operator $1 - \mc$ 
on the Lie algebra of $H$ (see also \cite{Ca} if ${\rm char}({\bf k})\neq 2$ and is good for $G$).
Moreover, it is shown by M. Costantini \cite{Co}, again for a spherical conjugacy class $C$, that
the decomposition of the coordinate ring of $C$ as a $G$-module is almost entirely
determined by the element $\mc$ (see \cite[Theorem 3.22]{Co} for the precise statement). When
$G$ is simple, the explicit description of $\mc$ for each spherical conjugacy class $C$ in $G$ is
given in \cite{CCC}, and the complete list of the
$\mc$'s, for $C$ spherical, is given in \cite[Corollary 4.2]{Ca}. 

In the second part of the paper, we study the element $\mc$ for an arbitrary
 conjugacy class $C$ 
in $G$. After examining some properties of $\Wc$, we show in Corollary~\ref{co-mc-maximal-length} that for each conjugacy class $C$ in
$G$, $\mc \in W$ is the unique maximal length element in its conjugacy class in $W$. In particular,
$\mc$ is an involution. When $C$ is spherical, the fact that $\mc$ is an involution is
also proved in \cite[Remark 4]{CCC} and \cite[Theorem 2.7]{Ca}.
Let
\begin{equation}\label{eq-M}
\calM = \left\{m \in W: \; m \; \mbox{is the unique maximal length element in} \; \calO_m\right\},
\end{equation}
where $\calO_m$ is the conjugacy class of $m$ in $W$. Then $\mc \in \calM$ for every conjugacy class $C$ in $G$. 
It is thus desirable to study the set $\calM$.

Assume that $G$ is simple. Using arguments from \cite{Ca}, we give in Theorem~\ref{th-main}
the complete list of 
elements in $\calM$. When  ${\rm char}({\bf k}) \neq 2$ and is good for $G$ (see
Remark~\ref{re-also} for more detail),
the list of elements in $\calM$
coincides with that in \cite[Corollary 4.2]{Ca} for the elements $\mc$ for $C$ spherical, and one thus has
\begin{eqnarray}\label{eq-calM-1}
\calM & =&\{\mc \in W: \; C \; \mbox{is a spherical conjugacy class in}\; G\} \\
\label{eq-calM-2} &=& 
\{\mc \in W: \; C \; \mbox{is a conjugacy class in}\; G\}.
\end{eqnarray}

As the last part of the paper, we consider in $\S$\ref{sec-An} the case of $G = \SL$. For any conjugacy class $C$ in $\SL$ and  
any involution $w \in W \cong S_{n+1}$, we
show in Theorem~\ref{th-main-SL} that
\[
C \cap (BwB) \neq \emptyset \hs \mbox{iff} \hs l_2(w) \leq r(C),
\]
where $l_2(w)$ is the number of distinct $2$-cycles in the cycle decomposition of $w$, and 
$r(C)= \min\{\rank(g-cI): c \in {\bf k}\}$
for any  $g \in C$.
Theorem~\ref{th-main-SL} is proved in $\S$\ref{subsec-proof-1} 
using  a criterion of Ellers-Gordeev \cite[Theorem 3.20]{EG2}. The same
criterion of Ellers-Gordeev also leads to a simple condition for a conjugacy class ${\mathcal O}$ in $W$ to
lie entirely in $\Wc$ for a given conjugacy class $C$ in $\SL$. See Theorem~\ref{th-O-WC}.

A direct consequence of Theorem~\ref{th-main-SL} is the following explicit description of the element $\mc$ for every
conjugacy class $C$ of $\SL$: for a positive integer $p$, let $\left[\frac{p}{2}\right]$ be the largest integer that is less than or equal to 
$\frac{p}{2}$. Let $m_0 = 1$, and for an integer
$1 \leq l \leq \left[\frac{n+1}{2}\right]$, let $m_l 
\in S_{n+1}$ be the involution with the cycle decomposition
\[
m_l = (1, n+1) (2, n) \cdots (l,\; n+2-l).
\]
Our Corollary~\ref{co-main-SL} says that 
for any conjugacy class $C$ in $SL(n+1, {\bf k})$,
\[
\mc =  m_{{}_{\left[ \frac{n+1}{2}\right]}}\hs \mbox{if} \; \;  r(C) \geq \left[ \frac{n+1}{2}\right], \hs \mbox{and} \hs 
\mc=m_{r({\scriptscriptstyle C})} \hs \mbox{if} \; \; r(C) < \left[ \frac{n+1}{2}\right].
\]
Note that $m_{{}_{\left[ \frac{n+1}{2}\right]}}$ is the maximal element in $S_{n+1}$.
The explicit description of $\mc$ for an arbitrary conjugacy class $C$ in other classical 
groups will be given in \cite{Ch}.

%{\bf FIX HERE: Using again the criterion by Ellers-Gordeev \cite{EG2}, we given in Corollary~\ref{co-O-C}
%As another consequence of Theorem~\ref{th-main-SL}, we show in Corollary~\ref{co-Ow-SL} that for any conjugacy class $C$ in $\SL$, 
%if $\Wc$ contains  an involution $w \in S_{n+1}$, then
%$\Wc$ contains the whole conjugacy class of $w$ in $S_{n+1}$. }

Assume that ${\rm char}({\bf k}) \neq 2$. 
Combining Theorem~\ref{th-main-SL} and \cite[Theorem 2.7]{Ca}, we show in
Corollary~\ref{co-spherical-SL} that $\Wc = \{w \in S_{n+1}:  w^2 = 1, l_2(w) \leq r(C)\}$ for a spherical conjugacy class $C$
in $\SL$.

\subsection{Notation}\label{subsec-notation}
Let $\Delta$ be the set of all roots of $G$ with respect to $H$, let $\Delta^+\subset \Delta$ 
be the set of positive roots 
determined by $B$, and let
$\Gamma$ be the set of simple roots in $\Delta^+$. We also write $\alpha > 0$ (resp. $\alpha < 0$)
if $\alpha \in \Delta^+$ (resp. $\alpha \in -\Delta^+$). 

For $\alpha \in \Gamma$, let $s_\alpha \in W$ be the reflection determined by $\alpha$.
The length function on $W$ is denoted by $l$. For a subset $W_1 \subset W$, an element $w \in W_1$
is called a minimal (resp. maximal) length element if $l(w) \leq l(w_1)$ (resp. $l(w) \geq l(w_1)$)
for all $w_1 \in W_1$.
For a subset $J$ of $\Gamma$, let
$W_J$ be the subgroup of $W$ generated by $\{s_\alpha: \alpha \in J\}$, let $w_{0, J}$ be the
maximal length element in $W_J$, and let $W^J \subset W$ be the set of minimal length representatives of $W/W_J$.

For $\delta \in {\rm Aut}(W)$, define the $\delta$-twisted conjugation of $W$
on itself  by $u \cdot_\delta v = \delta(u) v u^{-1}$ for $u, v \in W,$
and for $w \in W$, let $\calO^{\delta}_w=\{\delta(u) w u^{-1}: u \in W\}$ and call it the $\delta$-twisted conjugacy class of $w$ in $W$. 
Let $w_0=w_{0, \Gamma}$ be the longest element in $W$. We will use $\delta_0$ both for the permutation on 
$\Delta^+$ given by $\delta_0(\alpha) = -w_0(\alpha)$ and for the  
automorphism of $W$ given by 
$\delta_0:  W \longrightarrow W:  \delta_0(w) = w_0 w w_0$ for $w \in W$.

\subsection{Acknowledgments} We thank X.-H. He  and
E. Ellers for helpful discussions and for
drawing our attention to various references. We are especially grateful to G. Carnovale for
answering our questions on spherical conjugacy classes and for her input to
Remark~\ref{re-also}.
We thank the referees for helpful suggestions for the revised version of the paper.
Our research was partially supported by HKRGC grants 703405 and 703407.

\section{The sets $\Wc$ and $\Wcm$ and the elements $\mc$}\label{sec-mc}

\subsection{$\Wcm$ in terms of $\mc$}\label{subsec-2.1}
Let the notation be as in $\S$\ref{subsec-1.1} and $\S$\ref{subsec-notation}. In particular, for each conjugacy class $C$ in $G$, the
subsets $\Wc$ and $\Wcm$ of $W$ are given  in (\ref{eq-Wc}) and (\ref{eq-Wcm}).

\begin{lemma}\label{le-Wc-Wcm-1}
One has $\Wc \subset \Wcm$ for every conjugacy class $C$ in $G$.
\end{lemma}

\begin{proof}
Let $w \in W$. If $C \cap (BwB) \neq \emptyset$, then $C \cap (Bw) \neq \emptyset$, so
$C \cap (BwB^-) \neq \emptyset$.
\end{proof}

\begin{lemma}\label{le-deodhar}
For any $w \in W$, 
\[
BwB^-B = \bigsqcup_{w' \in W, w \leq w'} Bw'B.
\]
\end{lemma}

\begin{proof} Clearly $BwB^-B$ is the union of some $(B, B)$-double cosets. Let $w' \in W$.
Then 
\[
B w' B \subset BwB^-B  \;\;\Longleftrightarrow \;\;(Bw'B) \cap (BwB^-B) \neq \emptyset 
\;\;\Longleftrightarrow \;\; (Bw'B) \cap (BwB^-) \neq \emptyset,
\]
which, by \cite[Corollary 1.2]{De}, is equivalent to $w \leq w'$.
\end{proof}

\begin{lemma}\label{le-Wc-Wcm-2}
Let $C$ be  a conjugacy class in $G$ and let $w \in W$. Then $w \in \Wcm$ if and only if $w \leq w'$ for some $w' \in \Wc$.
\end{lemma}

\begin{proof} Since $C$ is conjugation invariant,
\[
C \cap (BwB^-) \neq \emptyset \;\;\Longleftrightarrow \;\; C \cap (BwB^-B) \nem,
\]
 which, by Lemma~\ref{le-deodhar}, 
is equivalent to $w \leq w'$ for some $w' \in \Wc$.
\end{proof}

\begin{lemma}\label{le-def-mc}\cite[$\S$1]{CCC}
Let $C$ be a conjugacy class in $G$. Then 

1) there is a unique $\mc \in W$ such that $C \cap (B\mc B)$ is dense in $C$;

2) $w \leq \mc$ for every $w \in \Wc$.
\end{lemma}

\begin{theorem}\label{th-Wcm-mc} For every conjugacy class $C$ in $G$, $\Wcm = \{w \in W: \, w \leq \mc\}$.
\end{theorem}

\begin{proof} Let $w \in W$. If $w \leq \mc$, then $w \in \Wcm$ by Lemma~\ref{le-Wc-Wcm-2}. Conversely, if 
$w \in \Wcm$, then again by Lemma~\ref{le-Wc-Wcm-2}, $w \leq w'$ for some $w' \in \Wc$. Since $w' \leq \mc$
by Lemma~\ref{le-def-mc}, one has $w \leq \mc$.
\end{proof}

\begin{lemma}\label{le-wc-closure}
If $C$ and $C'$ are two conjugacy classes in $G$ such that $C' \subset \overline{C}$, then 
$\mcp\leq \mc$. 
\end{lemma}

\begin{proof}
Since $\emptyset\neq C'\cap (B\mcp B)\subset \overline{C}=\overline{C\cap(B \mc B)}\subset \overline{B \mc B},$ 
one has $\mcp \leq \mc$.
\end{proof}

A conjugacy class $C$ in $G$ is said to be central if $C = \{z\}$ for some $z$ in the center of $G$.
It is clear that $\mc = 1$ if $C$ is central.

\begin{lemma}\label{le-central} A conjugacy class $C$ in $G$ is central if and only if $\mc = 1$.
\end{lemma}

\begin{proof}  Assume that $\mc=1$.
Then $C \subset B$ by Lemma~\ref{le-def-mc}. Conjugating $C$ by a representative of
$w_0$  in $N_G(H)$, one has $C \subset B^-$, so $C \subset H$. Let $z \in C$ be arbitrary.
Let $\alpha$ be any root.
If $zuz^{-1} \neq u$ for some $u \in U_\alpha$, the one-dimensional unipotent subgroup of $G$ defined by $\alpha$, then
$u^{-1}zu = u^{-1} z u z^{-1} z \in U_\alpha z$ is not contained in $H$, contradicting the fact
that $C \subset H$. Thus $zuz^{-1} = u$ for every $u \in U_\alpha$. Consequently $z$ is in the center of $G$ and
$C = \{z\}$.
\end{proof}
 
\subsection{Some properties of $\Wc$ and $\Wcm$}

We recall some  results from \cite{EG1, GKP}.

\begin{definition}\label{de-ascent-strong-conj} 1) \cite[Definition 3.1]{EG1} Let $w, w' \in W$. An ascent from 
$w$ to $w'$ is a sequence $\{\alpha_j\}_{1 \leq j \leq k}$ in  $\Gamma$ such that
\[
w' = s_{\alpha_k} \cdots s_{\alpha_2} s_{\alpha_1} w \,s_{\alpha_1} s_{\alpha_2} \cdots s_{\alpha_k}
\]
and $l(s_{\alpha_j} \cdots s_{\alpha_2} s_{\alpha_1} w \,s_{\alpha_1} s_{\alpha_2} \cdots s_{\alpha_j}) \geq
l(s_{\alpha_{j-1}} \cdots s_{\alpha_2} s_{\alpha_1} w \, s_{\alpha_1} s_{\alpha_2} \cdots s_{\alpha_{j-1}})$
for every $1 \leq j \leq k$. Write $w' \longleftarrow w$ if there is an ascent from $w$ to $w'$ or if $w' = w$.

2) \cite[$\S$2.9]{GKP} For $w, w', x \in W$, write $w \stackrel{x}{\sim} w'$ if $l(w) = l(w')$, $w' = x w x^{-1}$,  and either 
$l(w')=l(xw) + l(x)$ or $l(w') = l(x) + l(w x^{-1})$. Write $w \sim w'$ if there exist sequences of $\{x_j\}_{1 \leq j \leq k}$
and $\{w_j\}_{1 \leq j \leq k}$ in $W$
such that 
\[
w \; \stackrel{x_1}{\sim} \; w_1 \; \stackrel{x_2}{\sim} \; \cdots \; \stackrel{x_k}{\sim} \; w_k = w'.
\]
\end{definition}

\begin{proposition}\label{pr-GKP}\cite[$\S$2.9]{GKP} Let $\calO$ be any conjugacy class in $W$.

1) For any $w \in \calO$, 
there exists a maximal length element $w' \in \calO$ such that $w' \longleftarrow w$;

2) If $w'$ and $w''$ are two maximal length elements in $\calO$, then $w' \sim w''$.
\end{proposition}

\begin{proposition}\label{pr-Wc} 
Let $C$ be a conjugacy class in $G$, and let $w, w' \in W$.

1) If $w' \longleftarrow w$ and $w \in \Wc$, then 
$w' \in \Wc$.

2) If $w \sim w'$ and $w \in \Wc$, then $w' \in \Wc$.
\end{proposition}

\begin{proof} 
1) is just \cite[Proposition 3.4]{EG1}. To see 2), assume that $w \stackrel{x}{\sim} w'$ for some $x \in W$, so
$w' = x w x^{-1}$,  and either 
$l(w')=l(xw) + l(x)$ or $l(w') = l(x) + l(w x^{-1})$. Assume first that $l(w')=l(xw) + l(x)$. Then
$B w'B = BxwBx^{-1}B$, so
\[
C \cap (B w'B) = C \cap (BxwBx^{-1}B) \supset C \cap (x w B x^{-1})  \neq \emptyset.
\]
Thus $C \cap (Bw' B) \neq \emptyset$, and so $w' \in \Wc$. The case of $l(w') = l(x) + l(w x^{-1})$ is proved similarly.
\end{proof}

\begin{remark}\label{re-WWWW} We
refer to \cite{EG1, EG2} for a more detailed study of the set $\Wc$ and in particular for the
case of $G = SL(n, {\bf k})$. It is proved in 
\cite{Ca1} by G. Carnovale that, when ${\rm char}({\bf k})\neq 2$ and is good for $G$,
 a conjugacy class $C$ in $G$ is spherical if and only if 
$\Wc$ consists only of involutions. 
\end{remark}

Recall that for $w \in W$, $\calO_w$ denotes the conjugacy class of $w$ in $W$. 

\begin{corollary}\label{co-O-Wcm}
Let $C$ be any conjugacy class in $G$, and let $w \in W$. If  $w \in \Wc$, then $\calO_w \subset \Wcm$.
\end{corollary}

\begin{proof} Suppose that $w \in \Wc$. By Proposition~\ref{pr-GKP} and Proposition~\ref{pr-Wc}, every maximal length element in 
$\calO_w$ is in $\Wc$. Since  maximal elements in $\calO_w$ with respect to the
Bruhat order are also maximal length elements \cite[Corollary 4.5]{He}, Corollary~\ref{co-O-Wcm} now follows from
Theorem~\ref{th-Wcm-mc}.
\end{proof}

\subsection{Some properties of $\mc$}

\begin{corollary}\label{co-mc-maximal-length}
For every conjugacy class $C$ in $G$, $\mc$ is the unique maximal length element in $\calO_{\mc}$.
\end{corollary}

\begin{proof}
By Proposition~\ref{pr-GKP}, 
there exists a maximal length element $w' \in \calO_{\mc}$ such that $w' \longleftarrow \mc$.
By Proposition~\ref{pr-Wc}, $w' \in \Wc$, so $w' \leq \mc$ by Lemma~\ref{le-def-mc}. 
Since $l(w') \geq l(\mc)$, one has $w' = \mc$.  
Thus  $\mc$ is a maximal length element in $\calO_{\mc}$. If $w_1$ is any maximal length element in $\calO_{\mc}$, then
$w_1 \in \Wc$ by Proposition~\ref{pr-GKP} and Proposition~\ref{pr-Wc}, and thus $w_1 \leq \mc$
by Lemma~\ref{le-def-mc}. Since $l(w_1) = l(\mc)$,
one has $w_1 = \mc$. Thus $\mc$ is the only maximal length element in $\calO_{\mc}$.
\end{proof}

\begin{corollary}\label{co-minimal length}
For any conjugacy class $C$ in $G$, $w_0 \mc$ is the unique minimal length 
element in the $\delta_0$-twisted conjugacy class $\calO^{\delta_0}_{w_0\mc}$ (see notation in $\S$\ref{subsec-notation}).
\end{corollary}

\begin{proof} Corollary~\ref{co-minimal length} follows from the fact that,
for any $w \in W$, the bijection
\begin{equation}\label{eq-phi}
\phi: \; \; W \longrightarrow W: \; \; u \longmapsto w_0 u, \hs u \in W,
\end{equation}
maps maximal length elements in $\calO_w$ to 
minimal length elements in $\calO^{\delta_0}_{w_0w}$.
\end{proof}

\begin{remark}\label{re-C}
Let $\tilde{G}$ be the connected and simply connected cover of $G$  with covering map  $\pi: \tilde{G} \to G$, 
and let $Z = \pi^{-1}(e)$, where $e$ is the identity element of $G$.
Let $\tilde{A} = \pi^{-1}(A)$, where $A \in \{H, B, B^-\}$.
 Identify the Weyl group
for $\tilde{G}$ with $W$. For any conjugacy class $C$ in $G$, $\pi^{-1}(C)$ is a union of conjugacy 
classes in $\tilde{G}$. Since $Z \subset \tilde{H} = \tilde{B} \cap \tilde{B}^-$, it is easy to see that for any conjugacy 
classes $\tilde{C}$ in $\pi^{-1}(C)$, $W_{\tilde{\sc}} = \Wc$ and $W_{\tilde{\sc}}^- = \Wcm$, and in particular,
$\mc = m_{\tilde{\sc}}$. Thus the subset $\{\mc:  C  \mbox{ a conjugacy class in } G\}$ of $W$ depends only on the
isogeny class of $G$.
\end{remark}

\subsection{$\delta$-twisted conjugacy classes}\label{subsec-twisted}

Let $\delta$ be an automorphism of $G$ such that $\delta(B) = B$ and $\delta(H) = H$,
and let $G$ act on itself by
\[
g \cdot_\delta h = \delta(g) h g^{-1}, \hs g, h \in G.
\]
The $G$-orbits for this action will be called $\delta$-twisted conjugacy classes in $G$.
In the study of a
certain Poisson structure on the
de Concini-Procesi compactification $\overline{G}$ of $G$ when $G$ is complex and of adjoint type \cite{EL-cplx},
one needs to know when a $\delta$-twisted conjugacy class in $G$ intersects with a Bruhat cell $BwB^-$.

Given a $\delta$-twisted conjugacy class $C_\delta$ of $G$, let
\begin{eqnarray}\label{eq-Wcd}
\Wcd &= &\{w \in W: \; C_\delta \cap (BwB) \neq \emptyset\},\\
\label{eq-Wcmd}
\Wcdm &= &\{w \in W: \; C_\delta \cap (BwB^-) \neq \emptyset\}.
\end{eqnarray}
All the arguments in $\S$\ref{subsec-2.1} carry through when $C$ is replaced by $C_\delta$. In particular, $\Wcd \subset \Wcdm$, and if 
$\mcd$ denotes the unique element in $W$ such that $C_\delta \cap (B\mcd B)$ is dense in $C_\delta$, then
\[
\Wcdm = \{w \in W: \; w \leq \mcd\}.
\]

Since $\delta(H) = H$ and $\delta(B) = B$, $\delta$  acts on the set $\Gamma$ of simple roots
and  on $W$ such that $\delta(w) = \delta \circ w \circ \delta^{-1}: \Delta \to \Delta$. 
Minimal length elements in $\delta$-twisted conjugacy classes in $W$ have been studied in \cite{GKP, He}. In particular,
for $w, w' \in W$, write $w' \longleftarrow^{\delta} w$ if 
$w_0w \longrightarrow_{F} w_0w'$ in the notation of \cite[Definition 2.2 (a)]{GKP}, and 
write $w \sim^\delta w'$ if $w_0w \sim_{F} w_0w'$ in the notation of \cite[Definition 2.2 (b)]{GKP},
where $F = \delta^{-1}\delta_0$.
Using the fact that the map $\phi$ in (\ref{eq-phi}) maps maximal length elements in a $\delta$-twisted conjugacy class
${\calO}_w^\delta$ to minimal length elements in the $\delta_0\delta$-twisted conjugacy class ${\calO}_{w_0w}^{\delta_0\delta}$, one sees that 
\cite[Theorem 2.6]{GKP} is equivalent to the following generalization of Proposition~\ref{pr-GKP}:

\begin{proposition}\label{pr-GKP-twisted}
Let $\calO^\delta$ be a $\delta$-twisted conjugacy class in $W$. 

1) For any $w \in \calO^\delta$, there exists a maximal length element $w' \in \calO^\delta$ such that $w' \longleftarrow^\delta w$;

2) If $w'$ and $w''$ are two maximal length elements in $\calO^\delta$, then $w' \sim^\delta w''$.
\end{proposition}

Similarly, we have the following generalization of Proposition~\ref{pr-Wc}.

\begin{proposition}\label{pr-Wc-twisted}
Let $C_\delta$ be a $\delta$-twisted conjugacy class in $G$, and let $w, w' \in W$.

1) If $w' \longleftarrow^\delta w$ and $w \in \Wcd$, then 
$w' \in \Wcd$.

2) If $w \sim^\delta w'$ and $w \in \Wcd$, then $w' \in \Wcd$.
\end{proposition}

\begin{proof} 
For 1), it is enough to assume that  $w' = s_{\delta(\alpha)} w s_\alpha$ for
$\alpha \in \Gamma$ and $l(w') \geq l(w)$. Assume also that $w(\alpha) \neq \pm \delta(\alpha)$ because
otherwise $w' = w$ and we are done. Fix a representative $\dot{s}_\alpha$ of $s_\alpha$ in 
$N_G(H)$. If $w' > ws_\alpha$, choosing $g \in C_\delta \cap Bw$, one has 
\[
\delta(\dot{s}_\alpha) g \dot{s}_\alpha^{-1} 
\in C_\delta \cap (B s_{\delta(\alpha)} Bws_\alpha B) = C_\delta \cap (Bw' B).
\]
Assume that $w' < ws_\alpha$. If $w' < s_{\delta(\alpha)} w$, then
$s_{\delta(\alpha)} (w(\alpha)) < 0$, and since $w(\alpha) \neq \pm\delta(\alpha)$, one has 
$w(\alpha) < 0$, which implies that 
$w' < ws_\alpha < w$, a contradiction. Hence $w' > s_{\delta(\alpha)} w$. Choose $g \in wB$. Then
\[
\delta(\dot{s}_\alpha) g \dot{s}_\alpha^{-1} 
\in C_\delta \cap (B s_{\delta(\alpha)} wBs_\alpha B) = C_\delta \cap (Bw' B).
\]

For 2), assume, in the notation of \cite[Definition 2.2 (b)]{GKP}, that  
$w_0w \stackrel{x}{\sim}_{\delta^{-1}\delta_0} w_0w'$ for some $x \in W$. Then 
$w' = \delta(y) w y^{-1}$, where $y = \delta^{-1}\delta_0(x^{-1})$, such that either
$l(w')=l(\delta(y)w) + l(y)$ or $l(w') = l(y) + l(w y^{-1})$. The rest of the proof of 2) is similar to 
that of 2) of Proposition~\ref{pr-Wc}. 
\end{proof}

\begin{corollary}\label{co-twisted}
For any $\delta$-twisted conjugacy class $C_\delta$ in $G$, $\mcd$ is the unique maximal length element in its
$\delta$-twisted conjugacy class in $W$.
\end{corollary}
 
\begin{proof} Using Proposition~\ref{pr-GKP-twisted} and Proposition~\ref{pr-Wc-twisted},
the proof of  the statement is similar to that of Corollary~\ref{co-mc-maximal-length}. 
%If $\delta^2 = 1$ on $\Gamma$, then by \cite[Corollary 7.7]{He}, $\mcd$ and $\delta(\mcd)^{-1}$ %are in the same 
%$\delta$-twisted conjugacy class in $W$ and have the same length, so $\mcd = \delta(\mcd)^{-1}$.
\end{proof}

\section{Conjugacy classes of $W$ with unique maximal length elements}\label{sec-unique-maximal length}
\subsection{The set $\calM$}\label{subsec-calM-def}
By Corollary~\ref{co-mc-maximal-length}, $\mc \in \calM$ for every conjugacy class $C$ in $G$,
where  $\calM$ is the subset of $W$ given in (\ref{eq-M}).
It is thus desirable to have a precise description of elements in $\calM$.
Clearly $\calM$ is in
one-to-one correspondence with conjugacy classes in $W$ that have unique maximal length elements.

It is easy to see that if $G = G_1 \times G_2 \times \cdots \times G_k$ is a product of simple groups
$G_j$ with Weyl groups $W_j$ for $1 \leq j \leq k$,  then 
\begin{equation}\label{eq-calM-00}
\calM = \calM_{1} \times \calM_{2} \times \cdots \times \calM_{k},
\end{equation}
where for $1 \leq j \leq k$, $\calM_j \subset W_j$ is defined as in (\ref{eq-M}).
It is thus enough to determine $\calM$ for $G$ simple. This will be done in $\S$\ref{subsec-calM}.

\begin{lemma}\label{le-m-0}
If $m \in \calM$, then $m^2 = 1$.
\end{lemma}

\begin{proof} By \cite[Corollary 3.2.14]{GP}, $m^{-1} \in \calO_m$. Since $l(m) = l(m^{-1})$, one has $m = m^{-1}$.
\end{proof}

\subsection{The correspondence between  $\calM'$ and $\calJ'$}
Introduce
\[
\calM' = \{m \in W: \; m^2 = 1 \; \mbox{and} \; m \; \mbox{is a maximal length element in} \; \calO_m\}.
\]
By Lemma~\ref{le-m-0}, $\calM \subset \calM'$. We first determine $\calM'$.

It is well-known that elements in $\calM'$
correspond to special subsets
of the set $\Gamma$ of simple roots. Indeed, 
minimal or maximal length elements in conjugacy classes of involutions in $W$ have been studied
(see, for example, \cite{GP, He, PR1, PR2, R1, S} and especially 
\cite[Remark 3.2.13]{GP} for minimal length elements, \cite[Theorem 1.1]{PR2}
for maximal length elements, and \cite[Lemma 3.6]{He} for minimal length elements in twisted conjugacy classes).  
We summarize the results on $\calM'$ in the following Proposition~\ref{pr-m}, and we
give a proof of Proposition~\ref{pr-m} for completeness.  

%\begin{lemma}\label{le-m-00}\cite[Lemma 3.2]{S85}
%Let $m \in W$ be an involution. If $\alpha \in \Gamma$ is such that %$l(s_\alpha m s_\alpha) = l(m)$, then $s_\alpha m s_\alpha = m$.
%\end{lemma}

\begin{definition}\label{de-J1} A subset $J$ of $\Gamma$ is said to have Property (1) if $J$ is $\delta_0$-invariant and 
$-w_0(\alpha) = -w_{0, \sJ}(\alpha)$ for all $\alpha \in J$.
Let $\calJ'$ be the collection of all subsets $J$ of $\Gamma$ that have Property (1).
For $J \in \calJ'$, let
$\mJ = w_0 w_{0, \sJ}$. 
For $m \in \calM'$, let $J_m = \{\alpha \in \Gamma: \; m(\alpha) = \alpha\}\subset \Gamma.$
\end{definition}

\begin{lemma}\label{le-Mp-0} 
If $m \in \calM'$, then $m = w_0w_{0, \sJ_m}$, and
$J_m\in \calJ'$.
\end{lemma}

\begin{proof}
By \cite[Lemma 3.2]{S85}, $s_{\alpha}ms_{\alpha}=m$ for any $\alpha \in \Gamma$ with $s_{\alpha}m>m$. By \cite[Proposition 3.5]{S85}, $m=w_0w_{0,J_m}$ and $J_m \in \mathcal{J}'$.
\end{proof}

\begin{proposition}\label{pr-m} 1) The map $\psi: \calM' \to \calJ': m \mapsto J_m$ is bijective with inverse given by
$J \mapsto \mJ$ for $J \in \calJ'$.

2) For $J, K \in \calJ'$, $\mJ$ and $\mK$ are in the same conjugacy class in $W$ if and only if there
exists $w \in W$ with $\delta_0(w) = w$ such that $w(J) = K$.
\end{proposition}

\begin{proof} 1) Since $m = w_0 w_{0, \sJm}$ for every $m \in \calM'$, $\psi$ is injective. To show that
$\psi$ is surjective, let $J \in \calJ'$ and we will prove that $\mJ \in \calM'$. 
Since $J$ is $\delta_0$-invariant, $\mJ$ is an involution. Property (1)
implies that $s_\alpha \mJ s_\alpha = \mJ$ for every $\alpha \in J$, so
$w \mJ w^{-1} = \mJ$ for every $w \in W_J$. Thus, if $u = w \mJ w^{-1}$ is an element in $\calO_{\mJ}$,
we can assume that $w \in W^J$ (see notation in $\S$\ref{subsec-notation}). Then
\begin{eqnarray*}
l(u) &\leq &l(w) + l(\mJ w^{-1}) =  l(w) + l(w_0) - l(w_{0, \sJ} w^{-1})\\
& = &
l(w) + l(w_0) - l(w_{0, \sJ}) - l(w^{-1}) = l(\mJ).
\end{eqnarray*}
This shows that $\mJ$ is of maximal length  in $\calO_{\mJ}$, so $\mJ \in \calM'$. To show that 
$\psi(\mJ) = J$, note that 
$J \subset J_{\mJ} = \{\alpha \in \Gamma: \mJ(\alpha) = \alpha\}$. 
If $\alpha \in \Gamma\backslash J$, then $w_{0, \sJ}(\alpha) > 0$ so $\mJ(\alpha)
=w_0w_{0, \sJ}(\alpha) < 0$. It follows that
$J_{\mJ} \subset J$. 
Thus $J_{\mJ} = J$, and $\psi(\mJ) = J$. This shows that
$\psi$ is surjective and that its inverse is given by $\psi^{-1}(J) = \mJ$.

2) Assume that $J, K \in \calJ'$ are such that $\mJ$ and $\mK$ are conjugate in $W$. Since $w \mJ w^{-1} = \mJ$ 
for any $w \in W_J$, we
may assume that  $\mK = w \mJ w^{-1}$ for some $w \in W^J$. Let $\alpha \in J$. It follows from $\mK w = w \mJ$ and $w \in W^J$ that
\[
\mK w(\alpha) = w \mJ(\alpha) = w(\alpha) > 0.
\]
Thus $w(\alpha) \in [K]^+$, where $[K]^+$ denotes the set of positive roots in the linear span of $K$.
Denote similarly by $[J]^+$ the set of positive roots in the linear span of $J$. Then $w([J]^+)\subset[K]^+$.
Since both $\mJ$ and $\mK$ are maximal length elements in the same  conjugacy class in $W$,
$l(\mJ) =l(\mK)$.  Since 
\[
l(\mJ) = l(w_0) - |[J]^+| \hs \mbox{and} \hs l(\mK) = l(w_0) - |[K]^+|,
\]
one has
$|[J]^+| =|[K]^+|$. Here for a set $A$, $|A|$ denotes the cardinality of $A$.
Thus $w([J]^+)=[K]^+$. It follows that $w(J) = K$. Now $\mK = w \mJ w^{-1}$ implies that
$w_{0, \sK} = \delta_0(w) w_{0, \sJ} w^{-1}$, so
$\delta_0(w) = w_{0, \sK} w w_{0, \sJ} = w$.

Conversely, if $J, K \in \calJ'$ are such that $w(J) = K$ for some $w \in W$ with $\delta_0(w) = w$, then
$w_{0, \sK} = w w_{0, \sJ} w^{-1} = \delta_0(w) w_{0, \sJ} w^{-1}$, so $\mK = w \mJ w^{-1}$.
\end{proof}

%\begin{remark}\label{re-proof-pr-m}
%The proof of 2) of Proposition~\ref{pr-m} also shows that for $J, K \in \calJ'$ and $w \in W^J$, if  $\mK = w \mJ w^{-1}$,
%then $w(J) = K$ and $\delta_0(w) = w$.
%\end{remark}

\subsection{The correspondence between $\calM$ and $\calJ$}\label{subsec-calM}
We now turn to the set $\calM$.
For a subset $J$ of $\Gamma$, an  $\alpha \in J$ is said to be {\it isolated} 
if $\la \alpha, \alpha'\ra = 0$ for every $\alpha' \in J\backslash\{\alpha\}$.
The following Definition~\ref{de-P2} is inspired by \cite[Lemma 4.1]{Ca}.

\begin{definition}\label{de-P2} A subset $J$ of $\Gamma$ is said to have 
Property (2) if for every isolated $\alpha \in J$, 
there is no $\beta \in \Gamma\backslash\{\alpha\}$
with the 
following properties:

a) $\la \alpha, \alpha \ra = \la \beta, \beta\ra$ and $\la \beta, \alpha\ra \neq 0$;

b) $\la \beta, \alpha'\ra = 0$ for all $\alpha' \in J\backslash \{\alpha\}$;

c) $-w_0(\beta) =\beta$.
\end{definition}

\begin{lemma}\label{le-m-000} If $m \in \calM$, then $J_m$ has
Properties (1) and (2).
\end{lemma}

\begin{proof} Let $m \in \calM$. By Lemma~\ref{le-m-0} and Lemma~\ref{le-Mp-0}, $J_m$ has Property (1).
Suppose that $\alpha \in J_m$ is an isolated point and that there exists $\beta \in \Gamma\backslash \{\alpha\}$ with
properties a),  b) and c) in Definition~\ref{de-P2}.  Let $J_m' = J_m\backslash \{\alpha\}$. Since $\alpha \in J_m$ is isolated,
one has $w_0(\alpha) = -\alpha$, so 
$m = w_0 s_\alpha w_{0, \sJmp}= s_\alpha w_0 w_{0, \sJmp},$
and by b) and c),
$m(\beta) = s_\alpha w_0 w_{0, \sJmp}(\beta)= s_\alpha w_0(\beta) = -s_\alpha(\beta) <0$. Thus
$s_\beta m s_\beta = m s_{m(\beta)} s_\beta = m s_\alpha s_\beta s_\alpha s_\beta.$
By a), $s_\alpha s_\beta s_\alpha s_\beta=s_\beta s_\alpha$, so 
$s_\beta m s_\beta = m s_\beta s_\alpha$, and thus
$s_\alpha s_\beta m s_\beta s_\alpha = s_\alpha m s_\beta.$
Since $l(s_\alpha m s_\beta) \geq l(s_\alpha m) -1 = l(m)$, and since $m\in \calM$, $s_\alpha m s_\beta=m$. It follows from $m s_\alpha m = s_\alpha$ that $s_\alpha s_\beta = 1$
which is 
a contradiction.

\end{proof}

Let $\calJ$ be the collection  of all subsets $J$ of $\Gamma$ with Properties (1) and (2).
A $J \in \calJ$ is said to be non-trivial if $J$ is neither empty nor the whole of $\Gamma$.

Identify $\Gamma$ with the
Dynkin diagram of $G$ and a subset $J$ of $\Gamma$ with a sub-diagram of the Dynkin diagram.
The following description of $\calJ$ for $G$ simple is obtained in \cite[Corollary 4.2]{Ca}.
We include the list here for the convenience of the reader.

\begin{lemma}\label{le-J-list}
Assume that $G$ is simple and of rank $n \geq 2$.
The following is a complete list of non-trivial 
$J \in \calJ$ 
with points in $J$ painted black:

\begin{enumerate}
\item $A_n$: $J_l=\left\{ \alpha_{i} : l+1 \leq i \leq n-l \right\}$ for $1 \leq l \leq \left[\frac{n+1}{2}\right]-1$:

\[
\xy
(-15,4)*{\alpha_1}; (0,4)*{\alpha_{l}};(15,4)*{\alpha_{l+1}}; (30,4)*{\alpha_{l+2}};(45,4)*{\alpha_{n-l-1}};(60,4)*{\alpha_{n-l}};(75,4)*{\alpha_{n-l+1}};(90,4)*{\alpha_{n}};
(-15,0)*{\circ}="1";(0,0)*{\circ}="2";(15,0)*{\bullet}="3"; (30,0)*{\bullet}="4"; (45,0)*{\bullet}="5"; (60,0)*{\bullet}="6";  (75,0)*{\circ}="7";(90,0)*{\circ}="8";
{\ar@{--} "1"; "2"  }; 
{\ar@{-} "2"; "3"  }; 
{\ar@{-} "3"; "4"  }; 
{\ar@{--} "4"; "5"  };
{\ar@{-} "5"; "6"  };  
{\ar@{-} "6"; "7"  };
{\ar@{--} "7"; "8"  };  
\endxy
\]

\item $B_n$: $J_{1, l}=\left\{ \alpha_i : l \leq i \leq n \right\}$ for $2 \leq l \leq n$:

\[
\xy
(0,4)*{\alpha_1};(15,4)*{\alpha_{l-1}}; (30,4)*{\alpha_{l}};(45,4)*{\alpha_{l+1}};(60,4)*{\alpha_{n-1}};(75,4)*{\alpha_{n}};
(0,0)*{\circ}="1";(15,0)*{\circ}="2"; (30,0)*{\bullet}="3"; (45,0)*{\bullet}="4"; (60,0)*{\bullet}="5";  (75,0)*{\bullet}="6";
{\ar@{--} "1"; "2"  }; 
{\ar@{-} "2"; "3"  }; 
{\ar@{-} "3"; "4"  }; 
{\ar@{--} "4"; "5"  };
{\ar@2{->} "5"; "6"  };  
\endxy
\]

\noindent
$J_{2, l}=\left\{\alpha_1, \alpha_3, \ldots, \alpha_{2l-1} \right\}\cup \left\{ \alpha_i : 2l+1 \leq i \leq n \right\}$, for $ 1 \leq l \leq \frac{n}{2}-1 $:

\[
\xy
(0,4)*{\alpha_1};(15,4)*{\alpha_{2}}; (30,4)*{\alpha_{2l-1}};(45,4)*{\alpha_{2l}};(60,4)*{\alpha_{2l+1}};(75,4)*{\alpha_{2l+2}};(90,4)*{\alpha_{n-1}};(105,4)*{\alpha_{n}};
(0,0)*{\bullet}="1";(15,0)*{\circ}="2"; (30,0)*{\bullet}="3"; (45,0)*{\circ}="4"; (60,0)*{\bullet}="5";  (75,0)*{\bullet}="6";(90,0)*{\bullet}="7";(105,0)*{\bullet}="8";
{\ar@{-} "1"; "2"  }; 
{\ar@{--} "2"; "3"  }; 
{\ar@{-} "3"; "4"  }; 
{\ar@{-} "4"; "5"  };
{\ar@{-} "5"; "6"  };  
{\ar@{--} "6"; "7"  };
{\ar@2{->} "7"; "8"  };  
\endxy
\]
If $n=2m$,
$J_3=\left\{\alpha_1, \alpha_3, \ldots, \alpha_{2m-3}, \alpha_{2m-1} \right\}$:
\[
\xy
(0,4)*{\alpha_1};(15,4)*{\alpha_{2}}; (30,4)*{\alpha_{2m-3}};(45,4)*{\alpha_{2m-2}};(60,4)*{\alpha_{2m-1}};(75,4)*{\alpha_{2m}};
(0,0)*{\bullet}="1";(15,0)*{\circ}="2"; (30,0)*{\bullet}="3"; (45,0)*{\circ}="4"; (60,0)*{\bullet}="5";  (75,0)*{\circ}="6";
{\ar@{-} "1"; "2"  }; 
{\ar@{--} "2"; "3"  }; 
{\ar@{-} "3"; "4"  }; 
{\ar@{-} "4"; "5"  };
{\ar@2{->} "5"; "6"  };  
\endxy
\]
If $n=2m+1$,
$J_4=\left\{\alpha_1, \alpha_3, \ldots, \alpha_{2m-1}, \alpha_{2m+1} \right\}$:
\[
\xy
(0,4)*{\alpha_1};(15,4)*{\alpha_{2}}; (30,4)*{\alpha_{2m-1}};(45,4)*{\alpha_{2m}};(60,4)*{\alpha_{2m+1}};
(0,0)*{\bullet}="1";(15,0)*{\circ}="2"; (30,0)*{\bullet}="3"; (45,0)*{\circ}="4"; (60,0)*{\bullet}="5";  
{\ar@{-} "1"; "2"  }; 
{\ar@{--} "2"; "3"  }; 
{\ar@{-} "3"; "4"  }; 
{\ar@2{->} "4"; "5"  };  
\endxy
\]

\item $C_n$:
$J_{1, l}=\left\{ \alpha_i : l \leq i \leq n \right\}$ for $2 \leq l \leq n$:

\[
\xy
(0,4)*{\alpha_1};(15,4)*{\alpha_{l-1}}; (30,4)*{\alpha_{l}};(45,4)*{\alpha_{l+1}};(60,4)*{\alpha_{n-1}};(75,4)*{\alpha_{n}};
(0,0)*{\circ}="1";(15,0)*{\circ}="2"; (30,0)*{\bullet}="3"; (45,0)*{\bullet}="4"; (60,0)*{\bullet}="5";  (75,0)*{\bullet}="6";
{\ar@{--} "1"; "2"  }; 
{\ar@{-} "2"; "3"  }; 
{\ar@{-} "3"; "4"  }; 
{\ar@{--} "4"; "5"  };
{\ar@2{<-} "5"; "6"  };  
\endxy
\]

\noindent
$J_{2, l}=\left\{\alpha_1, \alpha_3, \ldots, \alpha_{2l-1} \right\}\cup \left\{ \alpha_i : 2l+1 \leq i \leq n \right\}$ 
for $ 1 \leq l \leq \frac{n}{2}-1 $:

\[
\xy
(0,4)*{\alpha_1};(15,4)*{\alpha_{2}}; (30,4)*{\alpha_{2l-1}};(45,4)*{\alpha_{2l}};(60,4)*{\alpha_{2l+1}};(75,4)*{\alpha_{2l+2}};(90,4)*{\alpha_{n-1}};(105,4)*{\alpha_{n}};
(0,0)*{\bullet}="1";(15,0)*{\circ}="2"; (30,0)*{\bullet}="3"; (45,0)*{\circ}="4"; (60,0)*{\bullet}="5";  (75,0)*{\bullet}="6";(90,0)*{\bullet}="7";(105,0)*{\bullet}="8";
{\ar@{-} "1"; "2"  }; 
{\ar@{--} "2"; "3"  }; 
{\ar@{-} "3"; "4"  }; 
{\ar@{-} "4"; "5"  };
{\ar@{-} "5"; "6"  };  
{\ar@{--} "6"; "7"  };
{\ar@2{<-} "7"; "8"  };  
\endxy
\]
If $n=2m$,
$J_3=\left\{\alpha_1, \alpha_3, \ldots, \alpha_{2m-3}, \alpha_{2m-1} \right\}$:
\[
\xy
(0,4)*{\alpha_1};(15,4)*{\alpha_{2}}; (30,4)*{\alpha_{2m-3}};(45,4)*{\alpha_{2m-2}};(60,4)*{\alpha_{2m-1}};(75,4)*{\alpha_{2m}};
(0,0)*{\bullet}="1";(15,0)*{\circ}="2"; (30,0)*{\bullet}="3"; (45,0)*{\circ}="4"; (60,0)*{\bullet}="5";  (75,0)*{\circ}="6";
{\ar@{-} "1"; "2"  }; 
{\ar@{--} "2"; "3"  }; 
{\ar@{-} "3"; "4"  }; 
{\ar@{-} "4"; "5"  };
{\ar@2{<-} "5"; "6"  };  
\endxy
\]
If $n=2m+1$, $J_4=\left\{\alpha_1, \alpha_3, \ldots, \alpha_{2m-1}, \alpha_{2m+1} \right\}$:
\[
\xy
(0,4)*{\alpha_1};(15,4)*{\alpha_{2}}; (30,4)*{\alpha_{2m-1}};(45,4)*{\alpha_{2m}};(60,4)*{\alpha_{2m+1}};
(0,0)*{\bullet}="1";(15,0)*{\circ}="2"; (30,0)*{\bullet}="3"; (45,0)*{\circ}="4"; (60,0)*{\bullet}="5";  
{\ar@{-} "1"; "2"  }; 
{\ar@{--} "2"; "3"  }; 
{\ar@{-} "3"; "4"  }; 
{\ar@2{<-} "4"; "5"  };  
\endxy
\]

\item $D_{2m}$: $J_{1, l}=\left\{ \alpha_i: 2l-1 \leq i \leq 2m \right\}$ for $ 2 \leq l \leq m$:

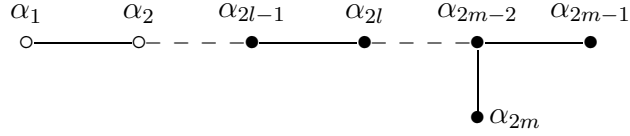
\[
\xy
(0,4)*{\alpha_1};(15,4)*{\alpha_{2}}; (30,4)*{\alpha_{2l-1}};(45,4)*{\alpha_{2l}};
(60,4)*{\alpha_{2m-2}};(75,4)*{\alpha_{2m-1}}; (65,-10)*{\alpha_{2m}};
(0,0)*{\circ}="1";(15,0)*{\circ}="2"; (30,0)*{\bullet}="3"; 
(45,0)*{\bullet}="4"; (60,0)*{\bullet}="5";  (75,0)*{\bullet}="6"; (60,-10)*{\bullet}="7";
{\ar@{-} "1"; "2"  }; 
{\ar@{--} "2"; "3"  }; 
{\ar@{-} "3"; "4"  }; 
{\ar@{--} "4"; "5"  };
{\ar@{-} "5"; "6"  };  
{\ar@{-} "5"; "7"  };  
\endxy
\]

\noindent
$J_{2, l}=\left\{\alpha_1, \alpha_3, \ldots, \alpha_{2l-1}\right\}\cup \left\{ \alpha_i: 2l+1 \leq i \leq 2m \right\}$
 for $ 1 \leq l \leq m-1$:

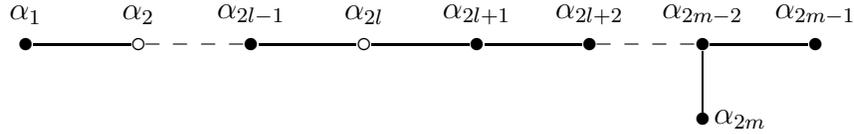
\[
\xy
(0,4)*{\alpha_1};(15,4)*{\alpha_{2}}; (30,4)*{\alpha_{2l-1}};(45,4)*{\alpha_{2l}};(60,4)*{\alpha_{2l+1}};
(75,4)*{\alpha_{2l+2}}; (90,4)*{\alpha_{2m-2}}; (105,4)*{\alpha_{2m-1}}; (95,-10)*{\alpha_{2m}};
(0,0)*{\bullet}="1";(15,0)*{\circ}="2"; (30,0)*{\bullet}="3"; (45,0)*{\circ}="4"; (60,0)*{\bullet}="5";  
(75,0)*{\bullet}="6"; (90,0)*{\bullet}="7"; (105,0)*{\bullet}="8";  (90,-10)*{\bullet}="9";
{\ar@{-} "1"; "2"  }; 
{\ar@{--} "2"; "3"  }; 
{\ar@{-} "3"; "4"  }; 
{\ar@{-} "4"; "5"  };
{\ar@{-} "5"; "6"  };  
{\ar@{--} "6"; "7"  };  
{\ar@{-} "7"; "8"  };
{\ar@{-} "7"; "9"  };     
\endxy
\]

\noindent
$J_3=\left\{\alpha_1, \alpha_3, \ldots, \alpha_{2m-3}, \alpha_{2m-1} \right\}$:

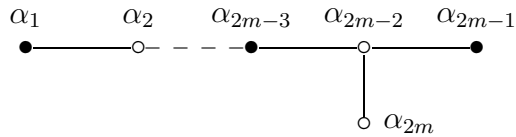
\[
\xy
(0,4)*{\alpha_1};(15,4)*{\alpha_{2}}; (30,4)*{\alpha_{2m-3}};(45,4)*{\alpha_{2m-2}};(60,4)*{\alpha_{2m-1}}; (51,-10)*{\alpha_{2m}};
(0,0)*{\bullet}="1";(15,0)*{\circ}="2"; (30,0)*{\bullet}="3"; (45,0)*{\circ}="4"; (60,0)*{\bullet}="5";  (45,-10)*{\circ}="6";
{\ar@{-} "1"; "2"  }; 
{\ar@{--} "2"; "3"  }; 
{\ar@{-} "3"; "4"  }; 
{\ar@{-} "4"; "5"  }; 
{\ar@{-} "4"; "6"  };
\endxy
\]

\noindent
$J_4=\left\{\alpha_1, \alpha_3, \ldots, \alpha_{2m-3}, \alpha_{2m} \right\}$:

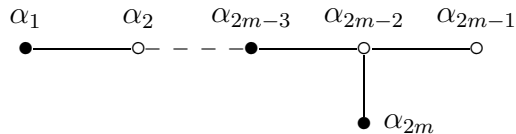
\[
\xy
(0,4)*{\alpha_1};(15,4)*{\alpha_{2}}; (30,4)*{\alpha_{2m-3}};(45,4)*{\alpha_{2m-2}};(60,4)*{\alpha_{2m-1}}; (51,-10)*{\alpha_{2m}};
(0,0)*{\bullet}="1";(15,0)*{\circ}="2"; (30,0)*{\bullet}="3"; (45,0)*{\circ}="4"; (60,0)*{\circ}="5";  (45,-10)*{\bullet}="6";
{\ar@{-} "1"; "2"  }; 
{\ar@{--} "2"; "3"  }; 
{\ar@{-} "3"; "4"  }; 
{\ar@{-} "4"; "5"  }; 
{\ar@{-} "4"; "6"  };
\endxy
\]

\item
$D_{2m+1}$: $J_{1, l}=\left\{ \alpha_i: 2l-1 \leq i \leq 2m+1 \right\}$ for $ 2 \leq l \leq m$:

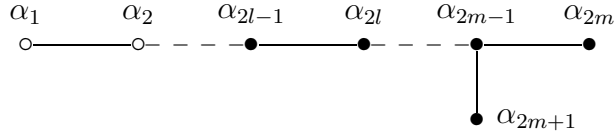
\[
\xy
(0,4)*{\alpha_1};(15,4)*{\alpha_{2}}; (30,4)*{\alpha_{2l-1}};(45,4)*{\alpha_{2l}};
(60,4)*{\alpha_{2m-1}};(75,4)*{\alpha_{2m}}; (68,-10)*{\alpha_{2m+1}};
(0,0)*{\circ}="1";(15,0)*{\circ}="2"; (30,0)*{\bullet}="3"; (45,0)*{\bullet}="4"; 
(60,0)*{\bullet}="5";  (75,0)*{\bullet}="6"; (60,-10)*{\bullet}="7";
{\ar@{-} "1"; "2"  }; 
{\ar@{--} "2"; "3"  }; 
{\ar@{-} "3"; "4"  }; 
{\ar@{--} "4"; "5"  };
{\ar@{-} "5"; "6"  };  
{\ar@{-} "5"; "7"  };  
\endxy
\]

\noindent
$J_{2, l}=\left\{ \alpha_1, \alpha_3, \ldots, \alpha_{2l-1} \right\} \cup 
\left\{ \alpha_i: 2l+1 \leq l \leq 2m+1 \right\}$ for $1 \leq l \leq m-1$:

\[
\xy
(0,4)*{\alpha_1};(15,4)*{\alpha_{2}}; (30,4)*{\alpha_{2l-1}};(45,4)*{\alpha_{2l}};
(60,4)*{\alpha_{2l+1}};(75,4)*{\alpha_{2l+2}}; (90,4)*{\alpha_{2m-1}}; (105,4)*{\alpha_{2m}}; (98,-10)*{\alpha_{2m+1}};
(0,0)*{\bullet}="1";(15,0)*{\circ}="2"; (30,0)*{\bullet}="3"; (45,0)*{\circ}="4"; 
(60,0)*{\bullet}="5";  (75,0)*{\bullet}="6"; (90,0)*{\bullet}="7"; (105,0)*{\bullet}="8";  (90,-10)*{\bullet}="9";
{\ar@{-} "1"; "2"  }; 
{\ar@{--} "2"; "3"  }; 
{\ar@{-} "3"; "4"  }; 
{\ar@{-} "4"; "5"  };
{\ar@{-} "5"; "6"  };  
{\ar@{--} "6"; "7"  };  
{\ar@{-} "7"; "8"  };
{\ar@{-} "7"; "9"  };     
\endxy
\]

\noindent
$J_3=\left\{\alpha_1, \alpha_3, \ldots, \alpha_{2m-1} \right\}$:

\[
\xy
(-15,4)*{\alpha_1};(0,4)*{\alpha_{2}};(15,4)*{\alpha_{2m-3}}; (30,4)*{\alpha_{2m-2}};(45,4)*{\alpha_{2m-1}};(60,4)*{\alpha_{2m}}; (52,-10)*{\alpha_{2m+1}};
(-15,0)*{\bullet}="1";(0,0)*{\circ}="2";(15,0)*{\bullet}="3"; (30,0)*{\circ}="4"; 
(45,0)*{\bullet}="5"; (60,0)*{\circ}="6";  (45,-10)*{\circ}="7";
{\ar@{-} "1"; "2"  }; 
{\ar@{--} "2"; "3"  }; 
{\ar@{-} "3"; "4"  }; 
{\ar@{-} "4"; "5"  }; 
{\ar@{-} "5"; "6"  };
{\ar@{-} "5"; "7"  };  

\endxy
\]

\item $E_6$:

\[
\xy
(-10,0)*{J_1:};(0,4)*{\alpha_1};(15,4)*{\alpha_3};(30,4)*{\alpha_4};(45,4)*{\alpha_5};(60,4)*{\alpha_6};(34,-10)*{\alpha_2};
(0,0)*{\bullet}="1";(15,0)*{\bullet}="3";(30,0)*{\bullet}="4";(45,0)*{\bullet}="5";(60,0)*{\bullet}="6";(30,-10)*{\circ}="2";
{\ar@{-} "1"; "3"  }; 
{\ar@{-} "2"; "4"  }; 
{\ar@{-} "3"; "4"  }; 
{\ar@{-} "4"; "5"  }; 
{\ar@{-} "5"; "6"  }; 
\endxy
\]

\[
\xy
(-10,0)*{J_2:};(0,4)*{\alpha_1};(15,4)*{\alpha_3};(30,4)*{\alpha_4};(45,4)*{\alpha_5};(60,4)*{\alpha_6};(34,-10)*{\alpha_2};
(0,0)*{\circ}="1";(15,0)*{\bullet}="3";(30,0)*{\bullet}="4";(45,0)*{\bullet}="5";(60,0)*{\circ}="6";(30,-10)*{\circ}="2";
{\ar@{-} "1"; "3"  }; 
{\ar@{-} "2"; "4"  }; 
{\ar@{-} "3"; "4"  }; 
{\ar@{-} "4"; "5"  }; 
{\ar@{-} "5"; "6"  }; 
\endxy
\]

\item $E_7$:

\[
\xy
(-10,0)*{J_1:};(0,4)*{\alpha_1};(15,4)*{\alpha_3};(30,4)*{\alpha_4};(45,4)*{\alpha_5};
(60,4)*{\alpha_6};(75,4)*{\alpha_7};(34,-10)*{\alpha_2};
(0,0)*{\circ}="1";(15,0)*{\bullet}="3";(30,0)*{\bullet}="4";(45,0)*{\bullet}="5";
(60,0)*{\bullet}="6";(75,0)*{\bullet}="7";(30,-10)*{\bullet}="2";
{\ar@{-} "1"; "3"  }; 
{\ar@{-} "2"; "4"  }; 
{\ar@{-} "3"; "4"  }; 
{\ar@{-} "4"; "5"  }; 
{\ar@{-} "5"; "6"  }; 
{\ar@{-} "6"; "7"  }; 
\endxy
\]

\[
\xy
(-10,0)*{J_2:};(0,4)*{\alpha_1};(15,4)*{\alpha_3};(30,4)*{\alpha_4};(45,4)*{\alpha_5};
(60,4)*{\alpha_6};(75,4)*{\alpha_7};(34,-10)*{\alpha_2};
(0,0)*{\circ}="1";(15,0)*{\bullet}="3";(30,0)*{\bullet}="4";(45,0)*{\bullet}="5";
(60,0)*{\circ}="6";(75,0)*{\bullet}="7";(30,-10)*{\bullet}="2";
{\ar@{-} "1"; "3"  }; 
{\ar@{-} "2"; "4"  }; 
{\ar@{-} "3"; "4"  }; 
{\ar@{-} "4"; "5"  }; 
{\ar@{-} "5"; "6"  }; 
{\ar@{-} "6"; "7"  }; 
\endxy
\]

\[
\xy
(-10,0)*{J_3:};(0,4)*{\alpha_1};(15,4)*{\alpha_3};(30,4)*{\alpha_4};(45,4)*{\alpha_5};
(60,4)*{\alpha_6};(75,4)*{\alpha_7};(34,-10)*{\alpha_2};
(0,0)*{\circ}="1";(15,0)*{\bullet}="3";(30,0)*{\bullet}="4";(45,0)*{\bullet}="5";
(60,0)*{\circ}="6";(75,0)*{\circ}="7";(30,-10)*{\bullet}="2";
{\ar@{-} "1"; "3"  }; 
{\ar@{-} "2"; "4"  }; 
{\ar@{-} "3"; "4"  }; 
{\ar@{-} "4"; "5"  }; 
{\ar@{-} "5"; "6"  }; 
{\ar@{-} "6"; "7"  }; 
\endxy
\]

\[
\xy
(-10,0)*{J_4:};(0,4)*{\alpha_1};(15,4)*{\alpha_3};(30,4)*{\alpha_4};(45,4)*{\alpha_5};
(60,4)*{\alpha_6};(75,4)*{\alpha_7};(34,-10)*{\alpha_2};
(0,0)*{\circ}="1";(15,0)*{\circ}="3";(30,0)*{\circ}="4";(45,0)*{\bullet}="5";
(60,0)*{\circ}="6";(75,0)*{\bullet}="7";(30,-10)*{\bullet}="2";
{\ar@{-} "1"; "3"  }; 
{\ar@{-} "2"; "4"  }; 
{\ar@{-} "3"; "4"  }; 
{\ar@{-} "4"; "5"  }; 
{\ar@{-} "5"; "6"  }; 
{\ar@{-} "6"; "7"  }; 
\endxy
\]

\item $E_8$:

\[
\xy
(-10,0)*{J_1:};(0,4)*{\alpha_1};(15,4)*{\alpha_3};(30,4)*{\alpha_4};(45,4)*{\alpha_5};
(60,4)*{\alpha_6};(75,4)*{\alpha_7};(90,4)*{\alpha_8};(34,-10)*{\alpha_2};
(0,0)*{\bullet}="1";(15,0)*{\bullet}="3";(30,0)*{\bullet}="4";(45,0)*{\bullet}="5";
(60,0)*{\bullet}="6";(75,0)*{\bullet}="7";(90,0)*{\circ}="8";(30,-10)*{\bullet}="2";
{\ar@{-} "1"; "3"  }; 
{\ar@{-} "2"; "4"  }; 
{\ar@{-} "3"; "4"  }; 
{\ar@{-} "4"; "5"  }; 
{\ar@{-} "5"; "6"  }; 
{\ar@{-} "6"; "7"  }; 
{\ar@{-} "7"; "8"  }; 
\endxy
\]

\[
\xy
(-10,0)*{J_2:};(0,4)*{\alpha_1};(15,4)*{\alpha_3};(30,4)*{\alpha_4};(45,4)*{\alpha_5};
(60,4)*{\alpha_6};(75,4)*{\alpha_7};(90,4)*{\alpha_8};(34,-10)*{\alpha_2};
(0,0)*{\circ}="1";(15,0)*{\bullet}="3";(30,0)*{\bullet}="4";(45,0)*{\bullet}="5";
(60,0)*{\bullet}="6";(75,0)*{\bullet}="7";(90,0)*{\circ}="8";(30,-10)*{\bullet}="2";
{\ar@{-} "1"; "3"  }; 
{\ar@{-} "2"; "4"  }; 
{\ar@{-} "3"; "4"  }; 
{\ar@{-} "4"; "5"  }; 
{\ar@{-} "5"; "6"  }; 
{\ar@{-} "6"; "7"  }; 
{\ar@{-} "7"; "8"  }; 
\endxy
\]

\[
\xy
(-10,0)*{J_3:};(0,4)*{\alpha_1};(15,4)*{\alpha_3};(30,4)*{\alpha_4};(45,4)*{\alpha_5};
(60,4)*{\alpha_6};(75,4)*{\alpha_7};(90,4)*{\alpha_8};(34,-10)*{\alpha_2};
(0,0)*{\circ}="1";(15,0)*{\bullet}="3";(30,0)*{\bullet}="4";(45,0)*{\bullet}="5";
(60,0)*{\circ}="6";(75,0)*{\circ}="7";(90,0)*{\circ}="8";(30,-10)*{\bullet}="2";
{\ar@{-} "1"; "3"  }; 
{\ar@{-} "2"; "4"  }; 
{\ar@{-} "3"; "4"  }; 
{\ar@{-} "4"; "5"  }; 
{\ar@{-} "5"; "6"  }; 
{\ar@{-} "6"; "7"  }; 
{\ar@{-} "7"; "8"  }; 
\endxy
\]

\item $F_4$:
\[
\xy
(-10,0)*{J_1:};(0,4)*{\alpha_1};(15,4)*{\alpha_2};(30,4)*{\alpha_3};(45,4)*{\alpha_4};
(0,0)*{\bullet}="1";(15,0)*{\bullet}="2";(30,0)*{\bullet}="3";(45,0)*{\circ}="4";
{\ar@{-} "1"; "2"  }; 
{\ar@2{->} "2"; "3"  }; 
{\ar@{-} "3"; "4"  }; 
\endxy
\]

\[
\xy
(-10,0)*{J_2:};(0,4)*{\alpha_1};(15,4)*{\alpha_2};(30,4)*{\alpha_3};(45,4)*{\alpha_4};
(0,0)*{\circ}="1";(15,0)*{\bullet}="2";(30,0)*{\bullet}="3";(45,0)*{\bullet}="4";
{\ar@{-} "1"; "2"  }; 
{\ar@2{->} "2"; "3"  }; 
{\ar@{-} "3"; "4"  }; 
\endxy
\]

\[
\xy
(-10,0)*{J_3:};(0,4)*{\alpha_1};(15,4)*{\alpha_2};(30,4)*{\alpha_3};(45,4)*{\alpha_4};
(0,0)*{\circ}="1";(15,0)*{\bullet}="2";(30,0)*{\bullet}="3";(45,0)*{\circ}="4";
{\ar@{-} "1"; "2"  }; 
{\ar@2{->} "2"; "3"  }; 
{\ar@{-} "3"; "4"  }; 
\endxy
\]

\item $G_2$:
\[
\xy
(-10,0)*{J_1:};(0,4)*{\alpha_1};(10,4)*{\alpha_2};(0,0)*{\circ}="1";(10,0)*{\bullet}="2";
{\ar@3{<-} "1"; "2"  };, \endxy, \hs \hs \xy (-10,0)*{J_2:};(0,4)*{\alpha_1};(10,4)*{\alpha_2};(0,0)*{\bullet}="1";(10,0)*{\circ}="2";
{\ar@3{<-} "1"; "2"  }.
\endxy
\]
\end{enumerate}\end{lemma}

\begin{theorem}\label{th-main} When  $G$ is simple, the map $\psi: \calM \to \calJ: \; m \mapsto J_m$ is a bijection
with inverse given by $\calJ \to \calM: \; J \mapsto \mJ$.
\end{theorem}

\begin{proof} It is clear that $\psi$ is injective. To show that $\psi$ is surjective, let $J \in \calJ$. 
We need to show that $\mJ \in \calM$. 
Let $m$ be any maximal length element in $\calO_{\mJ}$.
By Proposition~\ref{pr-m}, $m = \mK$, where  $K \in \calJ'$ and  there exists
$w \in W^J$ such that $w(J) = K$ and $\delta_0(w) = w$. 

By examining the list of all non-trivial $J$'s in $\calJ$ in Lemma~\ref{le-J-list} and the cases of
 $J = \emptyset$ or $J = \Gamma$, every
$J \in \calJ$ uniquely embeds into the Dynkin diagram with
Property (1) except in the cases of  $J_{1, m}, J_3, J_4$ for $D_{2m}$ and $J_4$ for
 $E_7$. In these cases, one can use results in 
\cite{PR1} to check directly that $\mJ \in \calM$.

\end{proof}

\begin{remark}\label{re-also} Let $G$ be simple. If ${\bf k} = {\mathbb C}$, then by  \cite[Remark 4.3]{Ca}, 
$\mJ$ for every $J \in \calJ$ is equal to 
$\mc$ for some
spherical conjugacy class $C$ in $G$. Consequently (\ref{eq-calM-1}) and (\ref{eq-calM-2}) hold.  
If ${\bf k}$ is algebraically closed 
with ${\rm char}({\bf k}) \neq 2$ and good for $G$,
spherical conjugacy classes
in $G$ are studied in \cite{Ca, Ca1} and classified in \cite{Ca2}. 
It was pointed out to us by G. Carnovale and by one of the
referees that \cite[Remark 4.3]{Ca} is still valid when ${\mathbb C}$ is replaced by such a field ${\bf k}$, so
(\ref{eq-calM-1}) and (\ref{eq-calM-2}) still hold. 
It follows from $\S$\ref{subsec-calM-def} that (\ref{eq-calM-1}) and (\ref{eq-calM-2}) also hold when $G$ is a product of 
simple groups.
\end{remark}

\begin{remark}\label{re-Perkins}
When $G$ is simple, the  minimal and maximal length elements as
well as their numbers in conjugacy classes of
involutions in 
$W$ are described  in  \cite{PR1}. For the classical groups, the descriptions in
\cite{PR1} are in terms of (signed) $2$-cycles for involutions in $W$  as opposed to sub-diagrams of the
Dynkin diagram.
\end{remark}

\begin{remark}\label{re-EG}
In \cite{EG1999}, E. Ellers and N. Gordeev considered the
intersections $C \cap (B^- w B)$, where $C$ is a non-central conjugacy class in a proper 
Chevalley group or finite twisted Chevalley group $G$ over an arbitrary field, and $w \in W$ 
is a  product of simple reflections in which
each simple reflection appears at most once.
In particular, 
\cite[Theorem 1]{EG1999} implies that such an intersection is always nonempty. In the 
notation of our paper, one has the following corollary  of \cite[Theorem 1]{EG1999}, which can also be checked directly case by case.
\end{remark}

\begin{proposition}\label{pr-EG1999}  Assume that $G$ is simple. Then
$w \leq \mc$ for every
non-central conjugacy class $C$ in $G$ and every Coxeter element $w \in W$.
\end{proposition}

\section{The case of $G = SL(n+1, {\bf k})$}\label{sec-An}

In this section, for an arbitrary conjugacy class $C$ in $\SL$, we give an explicit 
condition for $C \cap (BwB) \neq \emptyset$ when $w \in W \cong  S_{n+1}$ is an involution. In particular, 
we describe $\mc\in S_{n+1}$ explicitly for every $C$.

\subsection{Notation}\label{subsec-nota-SL}
As is standard, take the Borel subgroup $B$ (resp. $B^-$) to consist of
all
upper-triangular (resp. lower-triangular) matrices in $SL(n+1, {\bf k})$, so that
$H = B \cap B^-$ consists of all diagonal matrices in $SL(n+1, {\bf k})$. 

Identify the Weyl group $W$ of $SL(n+1, {\bf k})$ with the group $S_{n+1}$ of 
permutations on the set of integers between $1$ and $n+1$.
For $w \in S_{n+1}$, let
\[
l_2(w) = |\{i \in [1, n+1]: \; w(i) > i\}|.
\]

%Note that 
%\[
%l_2(w) = \rank(\rho(w)-I_{n+1}),
%\]
%where  $\rho: S_{n+1} \to GL(n+1, {\bf k})$ is the standard
%representatin of $S_{n+1}$ on ${\bf k}^{n+1}$.

Every conjugacy class $C$ in $\SL$ contains some $g$ of (upper-triangular) Jordan form.
We define the eigenvalues for $C$ to be the eigenvalues of such a $g \in C$ and similarly define
the number and sizes of the Jordan blocks of $C$ corresponding to an eigenvalue.
Let $I_{n+1} \in \SL$ be the identity matrix. For $g \in GL(n+1, {\bf k})$, define 
\begin{eqnarray*}
d(g) &=& \max \{\dim \ker (g - cI_{n+1}): \;c \in {\bf k}\}\\
r(g) &= &n+1-d(g)= \min \{\rank (g - cI_{n+1}): c \in {\bf k}\}.
\end{eqnarray*}
For a conjugacy class $C$ in $SL(n+1, {\bf k})$, define $r(C) = r(g)$ for any $g \in C$, and let
\[
l(C)   =\min\left\{r(C), \, \left[\frac{n+1}{2}\right]\right\}.
\]
%Note that if $C$ is a unipotent conjugacy class, then $d(C)$ is the number of Jordan blocks in any
%$g \in C$ that is of Jordan form. On the other hand, if $C$ is a semi-simple conjugacy class, then
%$d(C)$ is the largest among the multiplicities of the eigenvalues of any $g \in C$.
%Two elements in $\SL$ are in the same conjugacy class in $\SL$ if and only if they are
%in the same conjugacy class in $GL(n+1, {\bf k})$. This fact will be used throughout  this
%section.

\subsection{The main theorem and its consequences}\label{subsec-main-SL}

\begin{lemma}\label{le-SL}
Let $C$ be a conjugacy class in $SL(n+1, {\bf k})$ and let $w \in S_{n+1}$.  If
$C \cap (BwB) \neq \emptyset$, then 
$l_2(w) \leq r(C)$.
\end{lemma}

\begin{proof}
Assume that $C \cap (BwB) \neq \emptyset$. Let $g \in C \cap (BwB)$, and write $g = b_1 \dot{w} b_2$, where
$b_1, b_2 \in B$ and $\dot{w}$ is any representative of $w$ in the normalizer of $H$ in $G$. Then for any 
non-zero $c \in {\bf k}$,
\[
\rank(g-cI_{n+1}) = \rank(b_1 \dot{w} b_2 - cI_{n+1}) = \rank(\dot{w}-cb_1^{-1}b_2^{-1}).
\]
Let $l = l_2(w)$ and label the elements in the set $\{i \in [1, n+1]: w(i) > i\}$
as $i_1, \ldots, i_l$ such that $w(i_1) > \cdots > w(i_l)$.
It is easy to see that for any $b \in B$, the columns of the matrix $\dot{w} -b$ corresponding to $i_1,\ldots,i_{l}$ are linearly independent.
 Indeed, for $1 \leq k \leq l$, let $v_k$ be the $i_k$'th column of
$\dot{w}-b$. Suppose that $x_1 v_1  + \cdots +x_{l} v_{l} = 0$, where $x_k \in {\bf k}$  for each 
$1 \leq k \leq l$. By looking at the $w(i_1)$'th coordinate of  $x_1 v_1  + \cdots +x_{l} v_{l}$, one knows that
$x_1 = 0$. Repeating this argument, one sees that $x_k = 0$ for every $1 \leq k \leq l$.
Consequently, $\rank(\dot{w}-b) \geq l=l_2(w)$ for every $b \in B$, and thus
$\rank(g-cI_{n+1}) \geq l_2(w)$ for every non-zero $c \in {\bf k}$. It follows that $r(C) = r(g) \geq l_2(w)$.
\end{proof}

Note that for an involution   $w \in S_{n+1}$, $l_2(w)$ is the number of $2$-cycles 
in the cycle decomposition of $w$, and $l_2(w) \leq \left[\frac{n+1}{2}\right]$.

\begin{theorem}\label{th-main-SL} Let $C$ be a conjugacy class in $SL(n+1, {\bf k})$ and let $w \in S_{n+1}$ be an involution. Then
$C \cap (BwB) \neq \emptyset$ if and only if 
$l_2(w) \leq l(C)$.
\end{theorem}

Theorem~\ref{th-main-SL} will be proved in $\S$\ref{subsec-proof-1} using a result of Ellers-Gordeev \cite{EG2}.
We have been informed by E. Ellers that in a yet to be published paper,  E. Ellers and N. Gordeev have proved Theorem~\ref{th-main-SL} 
for the case of $w = w_0$, the longest element in $S_{n+1}$, but for the groups $GL(n+1, {\bf k})$ and $SL(n+1, {\bf k})$ for an 
arbitrary field ${\bf k}$.

We now give some corollaries of Theorem~\ref{th-main-SL}. We first determine the element $\mc$ for every conjugacy class $C$ in $\SL$. 

List the
simple roots as $\Gamma = \{\alpha_1, \alpha_2, \ldots, \alpha_n\}$ in the standard way.
Recall that $w_0$ is the longest element in $S_{n+1}$ and that for a subset $J$ of $\Gamma$, $w_{0, J}$ is the
longest element in the subgroup of $S_{n+1}$ generated by simple roots in $J$.
For an integer $0 \leq l \leq \left[\frac{n+1}{2}\right]$, let $m_l = w_0 w_{0, J_l}$, 
where $J_l = \emptyset$ for $l =\left[\frac{n+1}{2}\right]$ and 
$J_l = \{\alpha_{l+1}, \ldots, \alpha_{n-l}\}$ for
$0 \leq l \leq \left[\frac{n+1}{2}\right]-1$.
Thus,
$m_0 = 1$, 
and $m_l$ has the cycle decomposition 
\[
m_l = (1, n+1) (2, n) \cdots (l,\; n+2-l), \hs \hs \mbox{for}\;\;\;1 \leq l \leq \left[\frac{n+1}{2}\right].
\]
In particular, $m_l = w_0$ for $l = 
\left[\frac{n+1}{2}\right]$.
Note that for $0 \leq l_1, l_2 \leq \left[\frac{n+1}{2}\right]$, 
$m_{l_1} \leq m_{l_2}$ if and only if $l_1 \leq l_2$.

\begin{corollary}\label{co-main-SL}
For any conjugacy class $C$ in $SL(n+1, {\bf k})$,
$\mc = m_{l({\scriptscriptstyle C})}$, i.e.,
$\mc = w_0$ if $r(C) \geq \left[ \frac{n+1}{2}\right]$ and $ 
\mc =m_{r({\scriptscriptstyle C})}$ if $r(C) < \left[ \frac{n+1}{2}\right]$.
\end{corollary}

\begin{proof} 
Let $C$ be any conjugacy class in $\SL$. By Corollary~\ref{co-mc-maximal-length}, Lemma~\ref{le-m-000} and Lemma~\ref{le-J-list}, $\mc = m_l$ for some
$0 \leq l \leq \left[\frac{n+1}{2}\right]$. Since $C \cap (B m_l B) \neq \emptyset$, $l \leq l(C)$ 
by Theorem~\ref{th-main-SL}. Since $C \cap (B m_{l({\scriptscriptstyle C})} B) \neq \emptyset$
by  Theorem~\ref{th-main-SL}, one also has $l(C) \leq l$. Thus  $l =l(C)$.
\end{proof}

\begin{corollary}\label{co-C-C}
Let $C$ and $C'$ be two conjugacy classes in $\SL$ such that $C'$ is contained in the closure of $C$. 
Let 
$w \in S_{n+1}$ be  an involution.
If $w \in W_{{\scriptscriptstyle C'}}$, then $w \in \Wc$.
\end{corollary}

\begin{proof} It follows from the definition that $r(C') \leq r(C)$, so $l(C') \leq l(C)$. Corollary~\ref{co-C-C} now follows
directly from Theorem~\ref{th-main-SL}.
\end{proof}

\begin{corollary}\label{co-Ow-SL} Let ${\mathcal O}$ be a conjugacy class of involutions in $S_{n+1}$
and let $C$ be a conjugacy class in $\SL$.
If ${\mathcal O}\cap  \Wc \neq \emptyset$, then $\calO \subset \Wc$.
\end{corollary}

\begin{proof} Since $l_2(w) = l_2(w')$ for every $w, w' \in \calO$, Corollary~\ref{co-Ow-SL} follows directly from
Theorem~\ref{th-main-SL}.
\end{proof}

Consider spherical conjugacy classes in $\SL$. Assume that ${\rm char}({\bf k}) \neq 2$.
By \cite[Theorem 3.2]{Ca2}, a spherical conjugacy class in $SL(n+1, {\bf k})$ is either 
semi-simple or unipotent. Moreover, a
non-central semi-simple conjugacy class $C$ in $SL(n+1, {\bf k})$ is spherical if and only if it has exactly
two distinct eigenvalues, and in this case, 
$r(C)$ is equal to the 
smaller multiplicity of the two eigenvalues.
A
unipotent conjugacy class $C$ in $SL(n+1, {\bf k})$ is spherical if and only if all of its Jordan blocks
are of size at most $2$, and in this case, $r(C)$ is the number of size $2$ Jordan blocks for $C$.
In particular, $l(C) = r(C)$ for every spherical conjugacy class in $\SL$.

\begin{corollary}\label{co-spherical-SL} Let $G = SL(n+1, {\bf k})$ with 
${\rm char}({\bf k}) \neq 2$. Then $\Wc = \{w \in S_{n+1}:  w^2 = 1, \,\,  l_2(w) \leq r(C)\}$ for every spherical conjugacy class $C$ in $\SL$.
\end{corollary}

\begin{proof} Let $C$ be a spherical conjugacy class in $\SL$. 
Since ${\rm char}({\bf k}) \neq 2$, \cite[Theorem 2.7]{Ca} applies, and one knows that if $w \in \Wc$, then $w$ is an involution, and by Theorem~\ref{th-main-SL}, 
$l_2(w) \leq r(C)$. Conversely, if $w \in S_{n+1}$ is an involution with $l_2(w) \leq r(C)$, then
$w \in \Wc$ by Theorem~\ref{th-main-SL}.
\end{proof}

\subsection{Partitions associated to conjugacy classes}
Recall (see for example \cite[Page 705]{EG2}) that for an integer $p > 0$, a partition of $p$ is a non-increasing sequence 
$\lambda = (\lambda_1, \ldots, \lambda_s)$ of positive integers such that $\lambda_1 + \cdots+ \lambda_s = p$, and
$s$ is called the length of $\lambda$. The shape of a partition $\lambda = (\lambda_1, \ldots, \lambda_s)$ of $p$ 
consists of $s$ rows of  empty boxes left-aligned with $\lambda_j$ boxes on the $j$-th row for each $1 \leq j \leq s$.
The partition $\lambda^*$ of $p$ whose shape is obtained from switching the rows and columns of the shape of $\lambda$ is
called the dual of $\lambda$.
Let $\lambda = (\lambda_1, \ldots, \lambda_s)$ and $\mu = (\mu_1, \ldots, \mu_t)$ be two partitions of $p$. Define
$\lambda \leq \mu$ if 
$\sum_{j=1}^k \lambda_j \leq \sum_{j=1}^k \mu_j$ for every $1 \leq k \leq t$. 
One has (see \cite[Section I.1.11]{Mc})
$\lambda \leq \mu$ if and only if $\mu^* \leq\lambda^*$.

For integers $p > 0$ and $0 \leq l \leq \left[\frac{p}{2}\right]$, let $(2^l, 1^{p-2l})=(2, \ldots, 2, 1, \ldots, 1)$ be the partition of
$p$  with $2$ appearing exactly $l$ times.

\begin{lemma}\label{le-lambda}
Let $p  > 0$ be an integer and let $0 \leq l \leq \left[\frac{p}{2}\right]$. Then for any partition $\mu = (\mu_1, \ldots, \mu_s)$
of $p$, $(2^l, 1^{p-2l}) \leq \mu$ if and only if $s \leq p-l$.
\end{lemma}

\begin{proof}
One has 
$(2^l, 1^{p-2l}) \leq \mu$ if and only if $\mu^* \leq (2^l, 1^{p-2l})^*$, and the latter is equivalent to $s \leq p-l$.
\end{proof}

We now recall the partition $\tilde{\nu}^*(C)$ of $n+1$ associated to a conjugacy class $C$ in $\SL$  given in
\cite[Section 3.4]{EG2}: let $c_1, \ldots, c_k$ be the distinct eigenvalues of $C$, and for $1 \leq j \leq k$, let
$d_j>0$  and 
$n_{j,1} \geq \cdots \geq n_{j, d_j}>0$ be respectively the number and the sizes of Jordan blocks with eigenvalue $c_j$.
Arrange the list of the eigenvalues so that $d_1 \geq \cdots \geq d_k$. Then $d_1 = n+1-r(C)$. The partition $\tilde{\nu}^*(C)$ of $n+1$
is given by $\tilde{\nu}^*(C)=(\xi_1, \ldots, \xi_{d_1})$, where for $1 \leq t \leq d_1$,
$\xi_t = n_{1, t} + \cdots + n_{k, t}$ (if $d_j < d_1$ for some $2 \leq j \leq k$, 
we set $n_{j, s} = 0$ for all $d_j < s \leq d_1$).

For $w \in S_{n+1}$, let $\lambda(w)$ be the partition of $n+1$ formed by the lengths of the cycles in the cycle
decomposition of $w$.

We are now ready to prove Theorem~\ref{th-main-SL}.

\subsection{Proof of Theorem~\ref{th-main-SL}}\label{subsec-proof-1}
Let $C$ be a conjugacy class in $\SL$ and let $w \in S_{n+1}$ be an involution.
If $C \cap (BwB) \neq \emptyset$, then $l_2(w) \leq r(C)$ by Lemma~\ref{le-SL}, so $l_2(w) \leq l(C)$.
Conversely, assume that  
$l_2(w) \leq l(C)$, or, equivalently, $l_2(w) \leq r(C)$. We need 
to show that $C \cap (BwB) \neq \emptyset$.
By \cite[Theorem 3.2.9(a)]{GP}, there exists a minimal length element $w' \in \calO_w$ 
 and an ascent from $w'$ to $w$. Thus, in the notation of \cite{EG2}, there is a 
tree $\Gamma_w$ with $w'\in T(\Gamma_w)$. By \cite[Theorem 3.20]{EG2}, it is enough to show that
$\lambda(w') \leq \tilde{\nu}^*(C)$. Since $\lambda(w') = (2^{l_2(w)}, 1^{n+1-2l_2(w)})$, 
by Lemma~\ref{le-lambda}, $\lambda(w') \leq \tilde{\nu}^*(C)$ if and only if $d_1 \leq n+1-l_2(w)$ which
is equivalent to $l_2(w) \leq r(C)$.
This proves Theorem~\ref{th-main-SL}.

\subsection{A condition for $\calO \subset \Wc$}
For a conjugacy class $\calO$ in $S_{n+1}$, let $\lambda(\calO) = \lambda(w)$ for any $w \in \calO$.
Another consequence of \cite[Theorem 3.20]{EG2} is the following generalization of Corollary~\ref{co-Ow-SL}.

\begin{theorem}\label{th-O-WC}
Let $\calO$ be a conjugacy class in $S_{n+1}$ and $C$ a conjugacy class in $\SL$. Then $\calO \subset \Wc$ if and only if
$\lambda(\calO) \leq \tilde{\nu}^*(C)$.
\end{theorem}

\begin{proof} By Proposition~\ref{pr-GKP} and Proposition~\ref{pr-Wc}, $\calO \subset \Wc$ if and only if all minimal length elements
of $\calO$ are in $\Wc$. If $w \in \calO$ is a minimal length element, let $\Gamma_w$ be the 
tree consisting only of $w$. By \cite[Theorem 3.20]{EG2}, $w \in \Wc$ if and only if $\lambda(\calO) \leq \tilde{\nu}^*(C)$.
\end{proof}

%%%%%%%%%%%%%%%%%%%%%%%%%%%%%%%%%%%%%%%%%%%%%%%%%%%%%%%%%%%%%%%%%%%%%%%%%%%%%%%%%%%%%%%%

%%%%%%%%%%%%%%%%%%%%%%%%%
\end{document}